\documentclass{article}

\usepackage{amsmath}
\usepackage{amssymb}
\usepackage{vmargin}
\setpapersize{USletter}
\setmarginsrb{2cm}{2cm}{2cm}{2cm}{0cm}{0cm}{1cm}{0.5cm}
\usepackage{hyperref}
\usepackage{palatino}
\usepackage{graphicx}
\usepackage{tabularx}
\usepackage{algorithm,algorithmic}

\newcommand{\leavethisout}[1]{}
\newcommand{\myquote}[1]{``#1''}
\newcommand{\mydefn}[1]{\emph{#1}}

\begin{document}

\title{Mathematics For Industry: A Personal Perspective}
\author{John Stockie\\
  Department of Mathematics, Simon Fraser University}
\date{}

\maketitle


\section{Introduction}
\label{sec:intro}

\myquote{I am an industrial mathematician.}

When asked to identify my profession or academic field of study, this is
the most concise answer I can provide.  However, this seemingly
straightforward statement is commonly greeted by a blank stare or an
uncomfortable silence, regardless of whether I am speaking to a fellow
mathematician or a non-mathematician.  I usually follow up with the
clarification: \myquote{I am an applied mathematician who derives much
  of my inspiration from the study of industrial problems that I
  encounter through collaborations with companies.}  This dispels some
confusion, but unfortunately still leaves a great deal open to
interpretation owing to the vagueness of the words \mydefn{mathematics},
\mydefn{industry} and \mydefn{company}, each of which covers an
extremely broad range of scientific or socio-economic activity.

To those academics who actually work in the field of
\emph{industrial mathematics} (and whose ``perspective'' referred to in
the title is the focus of this article) this ambiguity is familiar and
untroubling.  However, for anyone less acquainted with the work of 
industrial mathematicians, some clarification is desirable especially for
anyone who might be considering entering the field.  This essay
therefore aims to shed light upon the nature of research being
done at the interface between mathematics and industry, paying
particular attention to the following questions:
\begin{itemize}
\item[] \emph{What is industrial mathematics?}  for which a proper
  answer depends sensitively on whether the mathematician in question is
  employed in industry or in academia, since these two sectors have
  demands, time scales and motives that are largely incommensurate.
  Because I am a university professor, this article naturally emphasizes
  the viewpoint of the academic mathematician.

\item[] \emph{Where is industrial mathematics?}  I presume that most
  readers will take no offense with the statement
  \myquote{mathematics is everywhere} which reflects the fact that
  math pervades the everyday operations of many organizations.  Be that
  as it may, it is certainly not the case that every organization making
  use of mathematics is also a source of novel, challenging (and
  ultimately publishable) mathematical problems for an academic
  mathematician.  The quest for interesting problems is really the
  \emph{raison d'\^{e}tre} for the industrial mathematician and so I
  will discuss how a wide range of non-academic organizations spanning
  many industry sectors can give rise to stimulating mathematics \dots\
  oftentimes in unexpected places.

\item[] \emph{How does one do industrial mathematics?} which has two
  facets: the first concerns the actual process of mathematical
  problem-solving (which I won't attempt to address here); and the
  second relates to the more fundamental question of how one goes about
  initiating an industrial collaboration in the first place.  The latter
  is commonly the most difficult stumbling block for anyone attempting
  to enter the field, owing at least in part to the disparities in
  culture and basic motivations that exist between academia and
  industry.

\item[] \emph{Why (or more precisely, what value is there in doing)
    industrial mathematics?}  There are obvious and
  well-doc\-u\-ment\-ed economic and strategic benefits to industry in
  deploying advanced mathematical solutions to difficult
  problems~\cite{macsinet-roadmap-2004, lery-etal-2012,
    smai-report-2002, nap-2012,
    siam-mii-2012}.  
  However, perhaps more relevant to many readers of this article is the
  complementary question of what benefits academic mathematicians can
  derive from industry collaborations, and how such projects can enrich
  their careers in terms of research, teaching, training and mentorship,
  and personal professional development.
\end{itemize}
I will attempt to answer these questions by means of several case
studies drawn from my own experience in tackling mathematical problems
from industry.  As a result, this account is necessarily a personal
perspective that is influenced to a large extent by my own research
interests bridging partial differential equations (mainly of parabolic
type), scientific computing (mainly finite volume methods) and fluid
dynamics.  I should confess at the outset that the contents of my own
\mydefn{mathematical toolbox} place significant limitations on the class
of problems and industrial applications that I am capable of making
meaningful contributions to.
I also make no claim to be the first person who has wrestled with the
questions stated above.  Indeed, there are many excellent accounts of
the pressing need for mathematics in industry and the important role
that mathematicians can play in solving real-world
problems~\cite{beauzamy-2002, campbell-1924, ockendon-1996,
  ockendon-2008, tayler-1990}.\ 
In particular, a paper by Bohun~\cite{bohun-2014} appeared recently
(while this article was under review) that addresses related issues
regarding the role and importance of the field of industrial
mathematics, and conveys the underlying message that the field is more
of an art than a well-defined scientific discipline.

\section{What is industrial mathematics?}
\label{sec:indmath}


Before attempting to answer this first question, it is most helpful to
introduce working definitions for \mydefn{industry} and \mydefn{mathematics} in
order to clarify how these terms are being used in the context of this
essay:
\paragraph*{Industry:} Most everyday uses of the term \mydefn{industry}
follow the dictionary definition that emphasizes either \myquote{the
  process of making products by using machinery and factories} or
\myquote{a group of businesses that provide a particular product or
  service}~\cite{mw-industry}.  Here I assume a much broader and more
inclusive definition that encompasses all \emph{end-users} of
mathematics, which includes not only public and private sector companies
but also federal/\-provincial/\-municipal government agencies,
hospitals, foundations, not-for-profits, etc.  By extension, I will use
the word \mydefn{company} to refer to any such non-academic
organization.  This inclusive definition is consistent with that
advocated in the OECD Global Science Forum's \emph{Report on Mathematics
  in Industry} as \myquote{any activity of economic or social value,
  including the service industry, regardless of whether it is in the
  public or private sector}~\cite{oecd-gsf-2008}.

\paragraph*{Mathematics:} I will take a similarly broad view of
\mydefn{mathematics} to refer to any branch of pure or applied
mathematics, statistics or computational science\footnote{Since
  numerical algorithms are simply the expression in computer code of a
  mathematical idea or procedure.} that can be employed to solve a
problem of interest to industry.  Moreover, I will refer primarily
(although not exclusively) to advanced mathematics that is in the
purview of academic mathematicians, and which relates to the the
underlying mathematical structure of a problem as well as the derivation
of exact or approximate solutions that are demonstrably (provably?)
correct.  This is in contrast with routine applications of well-known
mathematical techniques, which is more commonly the approach used by
engineers and other applied scientists.

\subsection{Mathematics for industry}
\label{sec:mfi}

With the above definitions in hand, the term \mydefn{industrial
  mathematics} (IM for short) refers to any mathematical treatment of
problems that arise from industrial applications.  Because this can
encompass such a huge spectrum of mathematics, IM is not so much a field
or sub-field of its own as it is a mathematical \emph{modus operandi}.

In his 2013 CAIMS-Mprime Industrial Mathematics Prize lecture, Brian Wetton
stated succinctly that \myquote{industrial mathematics is mathematics
  that industry is willing to pay for}~\cite{wetton-caimsprize-2013}.
I would like to propose an expanded definition of industrial mathematics
that includes three classes of activity:
\paragraph*{Mathematics IN industry (or MII)} which refers to
mathematics that is done by non-academic mathematicians who work as
employees of a company.  This is the perspective taken in the MII report
commissioned by SIAM~\cite{siam-mii-2012} and also in articles such
as~\cite{beauzamy-2002,campbell-1924}.
\paragraph*{Mathematics FOR industry (or MFI)} which is performed by
academic mathematicians within universities as part of a collaboration
with a company.  There is no necessity that the company is actually
funding the research, but rather that the development of the mathematics
is driven by needs of the partner organization and that there is a
two-way collaboration between academics and industrialists.  The
terminology MFI is consistent with Tayler~\cite{tayler-1990}.
\paragraph*{Mathematics INSPIRED BY industry (or MIBI)} by which I refer
to the mathematical analysis of problems that arise in an industrial
setting, but where all of the mathematical work is performed in an
academic setting that is largely isolated from the demands, pressures
and constraints of industry and industry-university collaborations.  For
example, work that begins as MII or MFI can still be a rich source of
more industrially-motivated problems that are not of direct interest to
the indunstry partner but nonetheless give rise to academic publications
or other activities that fall into the realm of MIBI.

\vspace{0.3cm}

\noindent 
The boundaries between these three classes are fuzzy and there is
significant overlap.  Nevertheless, this division provides a useful
high-level separation of industrial mathematics that may be represented
symbolically as
\begin{gather*}
  \text{IM} = \text{MII} \;\cup\; \text{MFI} \;\cup\; \text{MIBI},
\end{gather*}
by making use of the acronyms above and slightly abusing set notation.
This terminology is my own, and the reader should be aware that it
conflicts with common usage in certain segments of the mathematical
community.  For example, in the UK and many parts of Europe, the term
\mydefn{maths-in-industry}~\cite{ockendon-2008} is often used to refer
to the last two classes of activity ($\text{MFI} \;\cup\; \text{MIBI}$),
which is actually the complement of MII according to the SIAM
definition!  Furthermore, a \emph{Journal of Mathematics \underline{in}
  Industry} and a \emph{European Consortium for Mathematics
  \underline{in} Industry} (ECMI) have scope that effectively covers all
of IM.

Using the above definitions, it is now possible to clearly identify the
focus of the remainder of this article as the subset of IM that arises
in the context of direct collaborations between academic mathematicians
and industry and hence sits squarely within the second class of activity
called \emph{mathematics for industry} (MFI).

\subsection{The industrial mathematician} 
\label{sec:imathematician}

Some insight into the desirable characteristics in an industrial
mathematician is afforded by the following quotations:
\begin{itemize}
\item[] \myquote{Almost by definition, industrial research is
  interdisciplinary} (Ockendon, \cite{ockendon-2008}).
\item[] \myquote{[Real life mathematics] requires barbarians: people
    willing to fight, to conquer, to build, to understand, with no
    predetermined idea about which tool should be used} (Beauzamy,
  \cite{beauzamy-2002}).
\end{itemize}
The first author indicates that, unlike in university departments,
problems cannot be neatly categorized into academic disciplines; rather,
real industrial problems tend to be highly interdisciplinary in nature
so that their solution requires expertise that crosses disciplinary
boundaries.  Therefore, industrial mathematicians must be willing and
able to synthesize knowledge from other fields and formulate this
knowledge mathematically.

The second quotation refers to the need for a diverse mathematical
toolkit, since industrial problems are typically not clean textbook
examples in which a single mathematical technique or numerical algorithm
is sufficient to obtain a solution.  As a result, the industrial
mathematician must also be well-versed in more than one area of
mathematics and should be flexible in their approach to problem-solving,
displaying a willingness and ability to use whichever technique is most
appropriate.  For example, it is common for academic applied
mathematicians to build a career on becoming an expert in a certain
method (or class of methods) and then hunting for problems that they can
apply their method(s) to.  This approach of \myquote{choosing the
  problem to suit the method} is in stark contrast to the industrial
mathematician who is typically forced to \myquote{choose the method
  best suited to solving the problem.}

These qualities of breadth and agility, both within mathematics and
outside of the discipline, are essential for a successful industrial
mathematician.  As the case studies in the next section will
demonstrate, one of the primary rewards for investing the time and
effort required to develop these qualities is to provide access to a
wide variety of challenging and frequently novel mathematical problems.

\leavethisout{
\begin{itemize}
\item[] \myquote{Industry needs mathematicians of an especially broad
    type---men whose interests naturally extend beyond their special
    field, and who are flexible enough to cooperate with
    non-mathematicians.  These \emph{industrial mathematicians} must
    inspire confidence by their firm grasp of physical realities, by the
    relevance of their mathematics, and by the ability to present their
    results clearly and convincingly.} (Campbell, \cite{campbell-1924})
\item[] \myquote{[Mathematics] rarely solves industrial problems,
    indeed most industrial problems are never solved and compromise
    decisions have to be made.} (Tayler, \cite{tayler-1990})
\end{itemize}
}

\section{Case Studies}
\label{sec:case-studies}

The remaining three questions posed in the Introduction concern the
``where, how and why'' of industrial mathematics, and will be addressed
through examples.  In particular, I have chosen five case studies from
among my own past and current industrial projects to illustrate the
breadth of problems encountered in industrial settings as well as the
stimulating mathematical questions that can arise.  These examples
exhibit tremendous diversity in a number of aspects:
\begin{itemize}
\item \emph{Mathematical techniques:} ranging over simple algebra and
  trigonometry, partial differential equations, finite volume methods,
  integral transforms and Green's functions, Bayesian inversion,
  asymptotic analysis, and homogenization theory.
\item \emph{Scientific disciplines:} including engineering, image
  processing, electrochemistry, fluid mechanics, atmospheric science,
  and plant physiology.
\item \emph{Industrial partners:} representing industry sectors of
  manufacturing, defense, automotive, alternate energy, mining, and
  agriculture.  The companies also cover the entire size spectrum,
  ranging from small engineering consulting firms, to medium-sized
  R\&D-intensive companies, multi-national corporations, and non-profit
  network of industry associations.
\end{itemize}
For each case study, I will only provide a brief overview of the
mathematical problem and solution, instead emphasizing aspects of each
problem that address the questions posed in the Introduction.  Wherever
possible, I will provide references in which the interested reader can
find more mathematical details.  I will also describe how the
collaboration was initiated and how the project was funded, which in
most cases was through some form of joint industry-government R\&D
funding scheme.

\subsection{Robotic pipe welding: A consulting contract}
\label{sec:welding}

The first case study arose from a consulting project with a Vancouver
engineering firm and was initiated through a personal introduction to
one of the partner engineers at a neighbourhood social
event\footnote{This is an excellent example of how fruitful industrial
  connections can arise in the most unexpected situations!}.  I was
drawn into conversation with this engineer who had recently entered into
a contract with a pipe welding machinery manufacturer to write the
control software for a robotic welder.  He was struggling with how to
prescribe the motion of the welding torch for a general class of joints
involving pipes of different radii and intersection angles.

\begin{figure}[hbt]
  \small
  \centering
  \begin{tabular}{ccccc}
    \includegraphics[height=0.16\textheight]{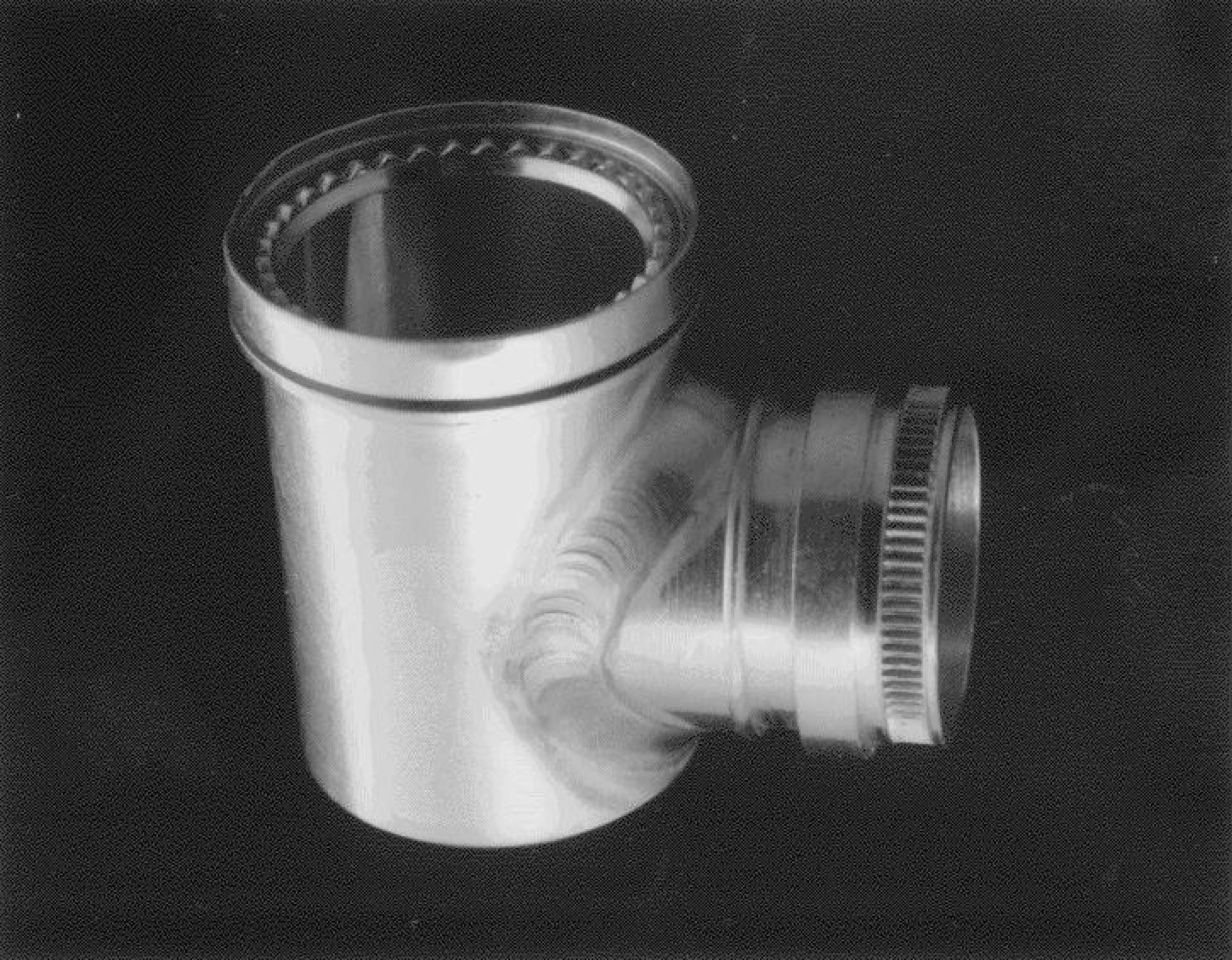}
    & \quad & 
    \includegraphics[height=0.18\textheight]{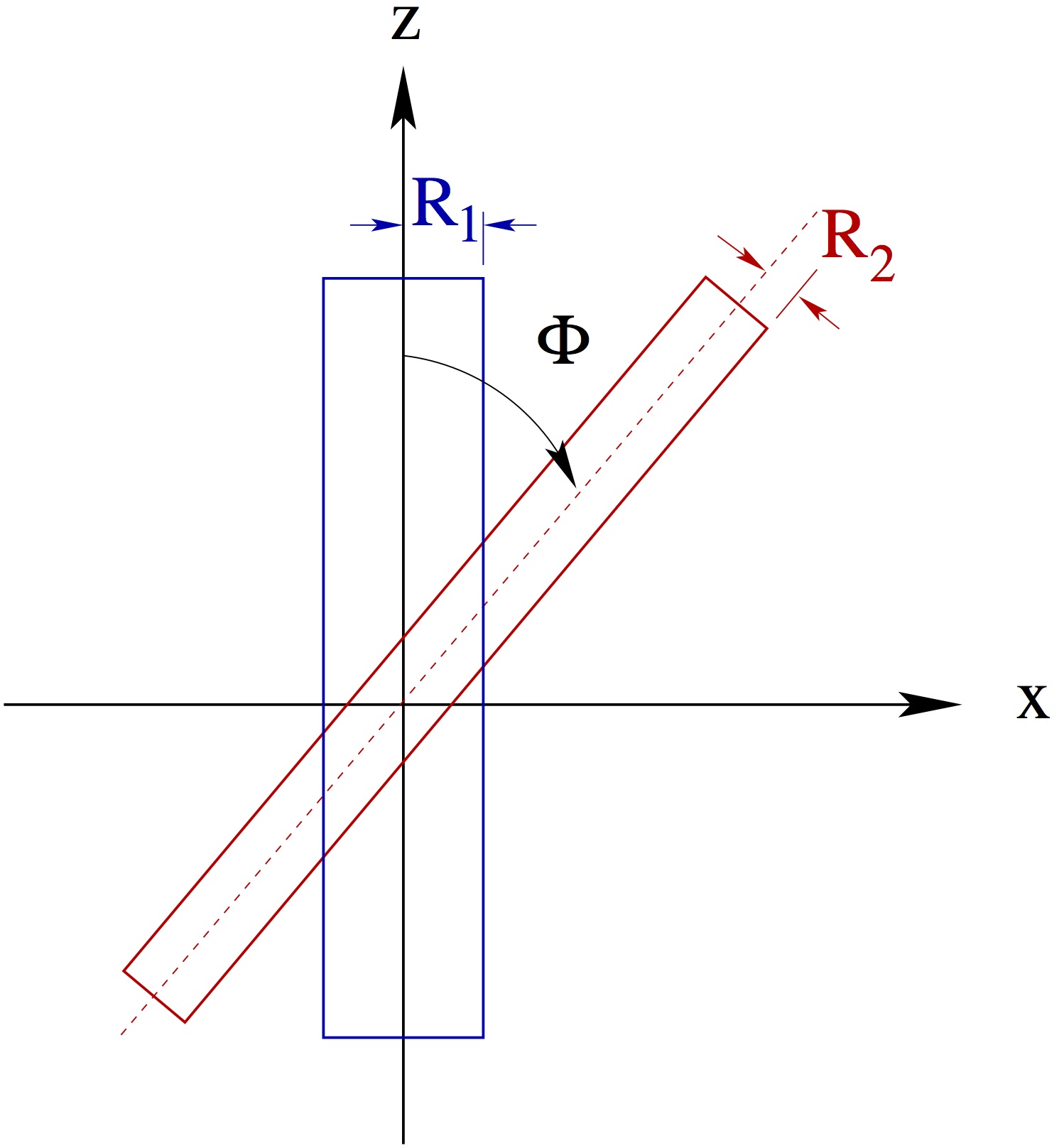}
    & & 
    \includegraphics[height=0.18\textheight]{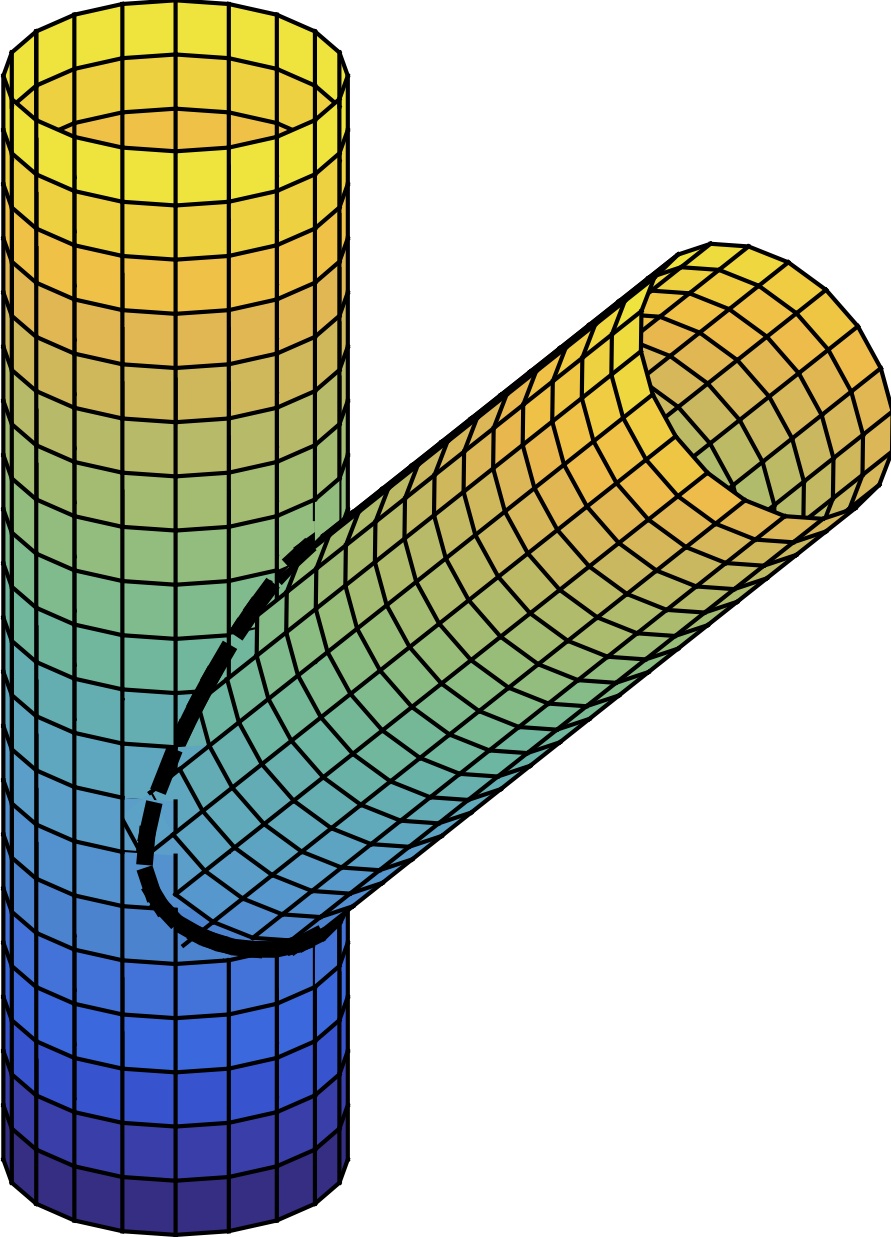}\\
    (a) && (b) && (c)
  \end{tabular}
  \caption{(a) Photograph of a pipe weld with joint angle
    $\Phi=90^{\circ}$. (b) Two cylinders with radii $R_1$ and $R_2$
    intersecting at an angle $\Phi$.  (c) The pipe weld solution is
    shown as a black line for values of $R_1=1$, $R_2=\frac{9}{10}$,
    $\Phi=45^\circ$.}
  \label{fig:pipe-weld}
\end{figure}

The problem geometry is illustrated in Figure~\ref{fig:pipe-weld}, which
includes a photograph of an actual pipe joint with a 90$^{\circ}$
intersection angle, as well as a diagram showing the general case where
two pipes with different radii ($R_1$ and $R_2$) are joined at an
arbitrary joint angle ($\Phi$).
I realized immediately that this was a straightforward exercise in
algebra and trigonometry that could be solved analytically and hence
easily implemented in software.  By writing a parametric description of
the pipes in cylindrical coordinates $(\theta_i,z_i)$ as follows
\begin{gather*}
  \text{\underline{Pipe 1:}}\;\;
    \begin{array}[t]{l}
    x = R_1 \cos\theta_1\\
    y = R_1 \sin\theta_1\\
    z = z_1
  \end{array}
  \qquad\qquad
  \text{\underline{Pipe 2:}}\;\;
  \begin{array}[t]{l}
    x = R_2 \cos\theta_2 \cos\Phi - z_2\sin\Phi\\
    y = R_2 \sin\theta_2\\
    z = R_2 \cos\theta_2 \sin\Phi + z_2\cos\Phi
  \end{array}
\end{gather*}
and equating $x,y$ and $z$ components, the two intersection
curves can be obtained in parametric form as
\begin{gather*}
  x = \pm \sqrt{R_1^2 - R_2^2 + R_2^2\cos^2\theta_2}, \qquad
  y = R_2\sin\theta_2, \qquad
  z = \frac{R_2\cos\theta_2 \mp \cos\Phi\sqrt{R_1^2 - R_2^2 +
      R_2^2\cos^2\theta_2}}{\sin\Phi} .
\end{gather*}
The solution derivation was a straightforward application of the Maple
software package\footnote{An expanded version of the Maple code was
  later posted on MapleSoft's Application Center at
  \url{http://www.maplesoft.com/applications/view.aspx?SID=3773}.}  and
formed the basis of an article appearing in the \emph{MapleTech}
journal~\cite{stockie-1998}. Some interesting implementation-related
issues arose surrounding automatic selection of the correct solution
branch as well as preventing the welding torch assembly from interfering
with the pipes at small joint angles.

This is the first and so far only consulting contract that I have been
involved with.  While the financial reward was attractive at the time
(especially for a PhD student living on a student-sized salary) this
reward has since been far surpassed by other more academic spin-offs.
First of all, I have used this problem to great effect in my
undergraduate teaching and outreach to high school students as a prime
example of the potential value to industry of expertise in mathematics.
Furthermore, I have subsequently been approached on a regular basis by
other welding companies or engineering firms to request assistance in
implementing the analytical pipe weld solution, when they encountered my
web page or the MapleSoft Application site.

\subsection{Land mine trip-wire detection: A study group problem}
\label{sec:tripwire}

\begin{figure}[hbt]
  \small
  \centering
  \begin{tabular}{cccc}
    \includegraphics[width=0.23\textwidth]{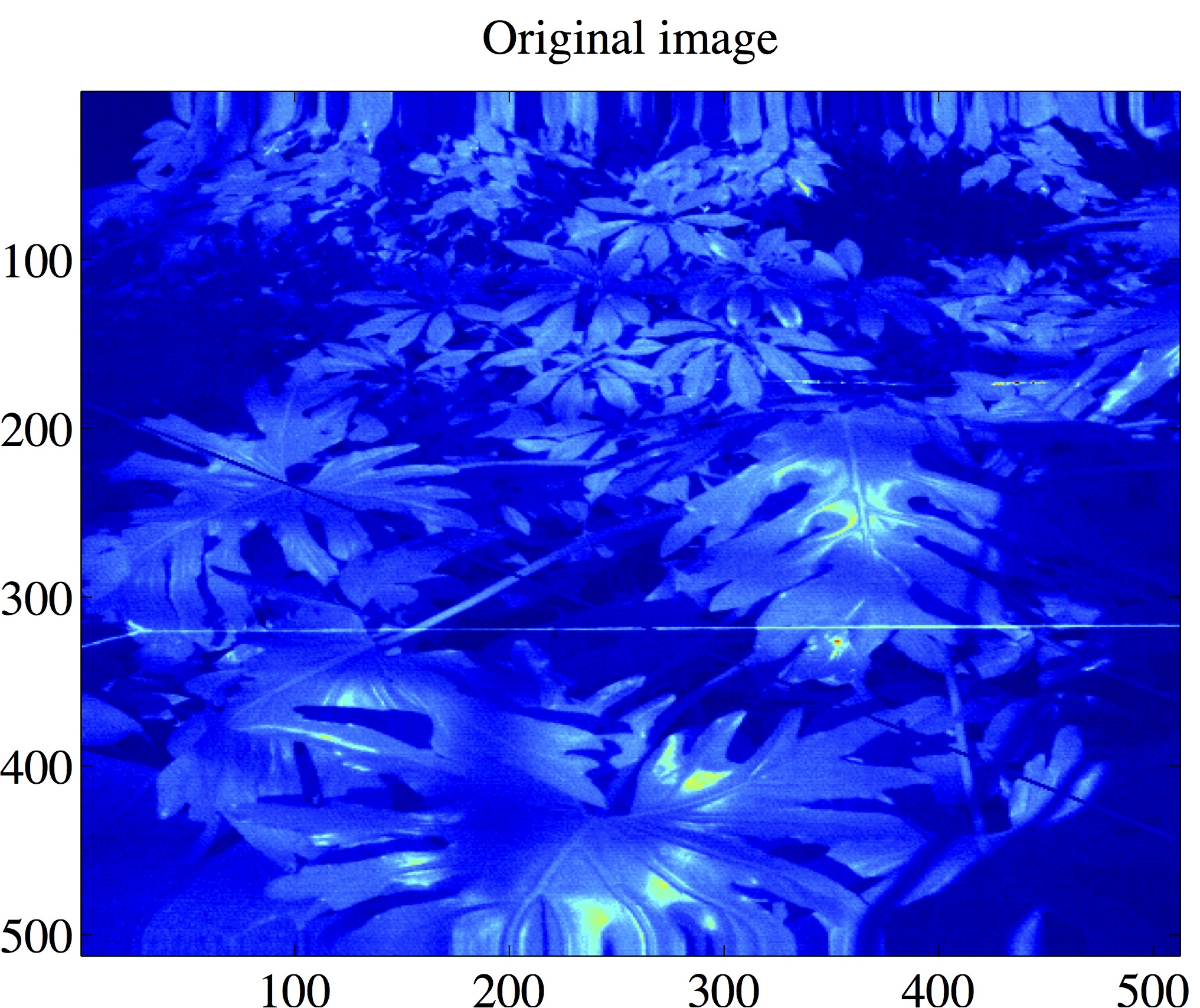} & 
    \includegraphics[width=0.23\textwidth]{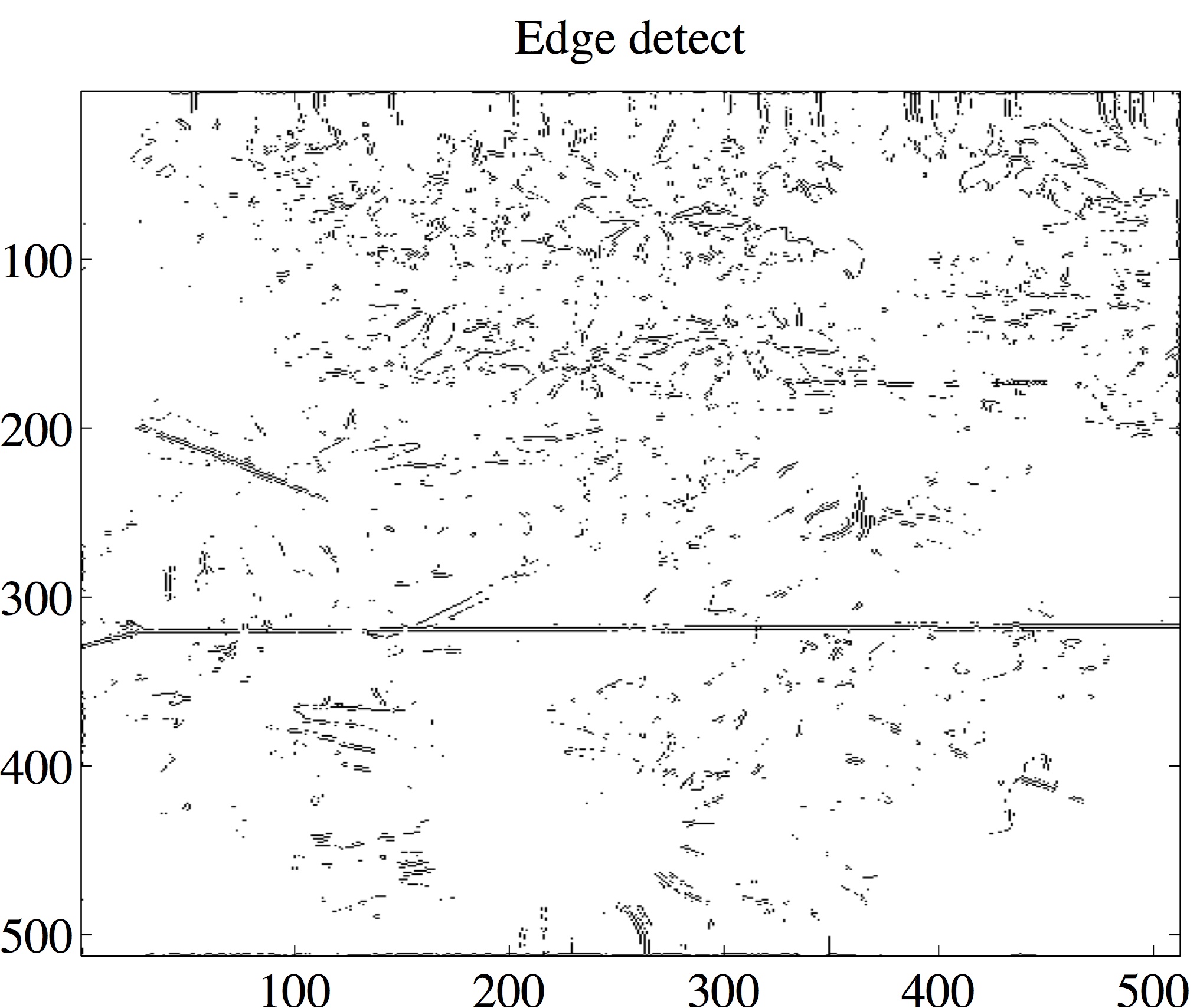} & 
    \includegraphics[width=0.24\textwidth]{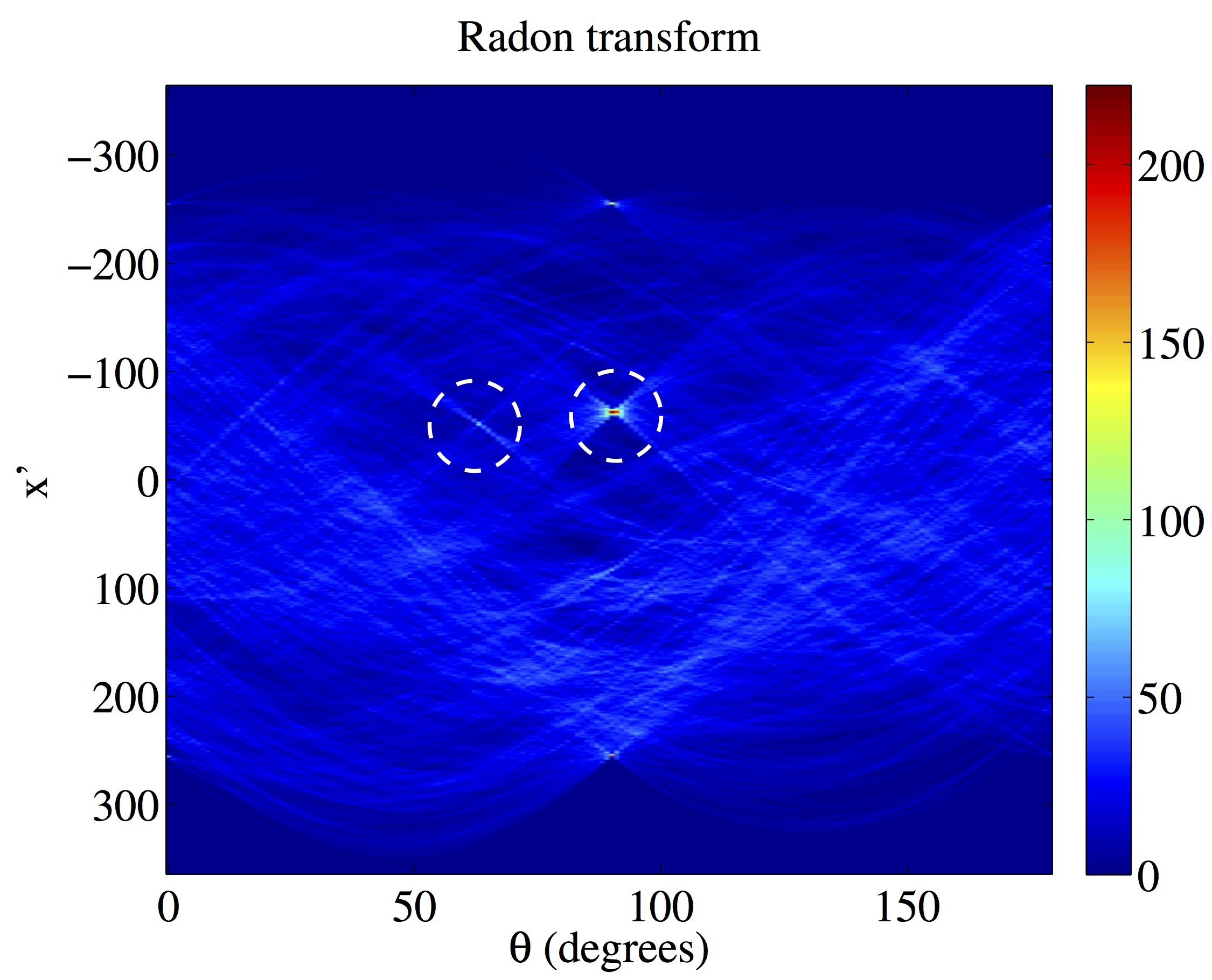} & 
    \includegraphics[width=0.23\textwidth]{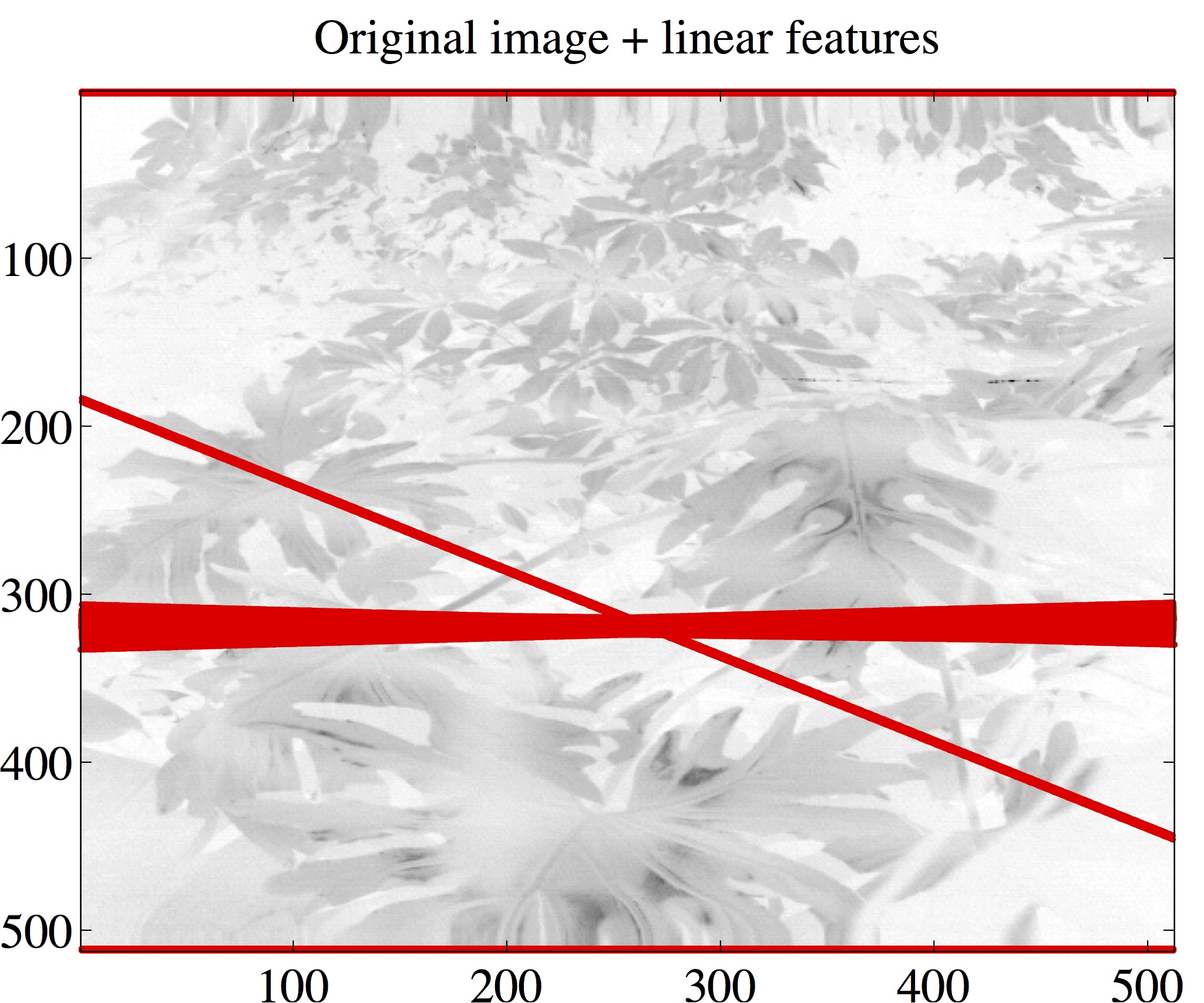}\\
    (a) & (b) & (c) & (d)
  \end{tabular}
  \caption{(a)~Original image containing two linear features: a bright
    horizontal line and a less obvious dark oblique line running
    downward from left to right.  (b)~Preprocessed image after Laplacian
    and edge detection filters.  (c)~Radon transform image $U(r,\theta)$
    with white circles around peaks that exceed the threshold $T$.
    (d)~Linear features (corresponding to the two peaks) are overlaid on
    the original image.}
  \label{fig:tripwire}
\end{figure}

The next case study is an image processing problem that was brought to
the second Industrial Problem-Solving Workshop (commonly known as a
study group outside of Canada) held in Calgary in 1998 by the Pacific
Institute for the Mathematical Sciences (PIMS).  The problem posed by
the industrial partner ITRES was to automatically identify land-mine
trip-wires within an image that is cluttered by vegetation, terrain,
man-made structures, etc.  A sample image file is provided in
Figure~\ref{fig:tripwire}(a).

Because trip-wires appear as linear features in an image, the team of
academics tackling this problem were naturally led to consider the \emph{Radon
  transform} that is a type of integral transform in which the integrals
are taken over straight lines.  If the intensity of a 2D image ${\cal
  I}$ is represented as a function $u(x,y)$ for $(x,y)\in {\cal I}$,
then the Radon transform may be written as
\begin{gather*}
  U(\rho,\theta) = \int_{\cal I} u(x,y) \, \delta( x\cos\theta -
  y\sin\theta - \rho)\, dx\, dy,
\end{gather*}
where $(\rho,\theta)$ are polar coordinates and $\delta(\cdot)$ is the
Dirac delta function.  Because the equation $x\cos\theta - y\sin\theta -
\rho=0$ represents a straight line in the original image, the above
integral has the effect that a linear feature with high intensity $u$
generates a single point $(\rho,\theta)$ at which the transformed
intensity $U$ is correspondingly large.  Hence, the largest values of
the Radon transform correspond to the primary linear features in the
original image.

When applying the Radon transform by itself to the test images provided,
it was not possible to reliably detect trip-wires.  Instead, it was
necessary to introduce a pre-processing step in which two filters
accentuate regions where the intensity changes rapidly.  The resulting
algorithm is summarized as Algorithm~\ref{alg:radon}.
\begin{algorithm}[hbt]
  \caption{Detecting linear features in an image using the Radon
    transform.}
  \label{alg:radon}
  \begin{algorithmic}[1]
    \STATE Apply \label{step:laplacian} a \emph{Laplacian filter}
    to each point in the raw image $u$ to accentuate regions of high
    curvature.
    \STATE Perform \label{step:edge} \emph{edge detection} on the
    filtered image.
    \STATE Calculate \label{step:radon} the \emph{Radon transform} $U$
    of the edge-detected image.
    \STATE Find \label{step:threshold} all points $(\rho,\theta)$ that
    \emph{exceed a threshold}, $U(\rho,\theta) > T$.
    \STATE Apply the \emph{inverse Radon transform} of each point found
    in step~\ref{step:threshold} to determine the linear features of
    interest.
  \end{algorithmic}
\end{algorithm}
The images obtained after steps~\ref{step:edge} and \ref{step:radon} of
the algorithm are shown in Figure~\ref{fig:tripwire}(b,c) and the final
linear features identified by the algorithm are overlaid on top of the
original image in Figure~\ref{fig:tripwire}(d).

These results were presented at the conclusion of the
Workshop~\cite{budd-stockie-1998} after which we had no further
communication with the industrial partner.  This is perhaps not
unexpected for a problem having such obvious military applications.
However, a publication with co-authors from ITRES did appear a few years
later in a conference proceedings~\cite{babey-etal-2000}, suggesting
that our results were eventually exploited by the company.  This is
another problem that has been highlighted in many outreach activities as
an example of the potential value of mathematics in
industry\footnote{Two examples of such outreach activities are at
  \url{http://plus.maths.org/content/saving-lives-mathematics-tomography}
  and \url{http://www.whydomath.org/node/tomography/landmine-wires.html}.}.

\subsection{Mathematical modelling of PEM fuel cells: A 
  multi-university collaboration} 
\label{sec:fuel-cells}

The next case study derives from a long collaboration with Ballard Power
Systems, an R\&D-intensive company based in Burnaby, British Columbia
that is a world leader in the development of polymer electrolyte
membrane (PEM) fuel cell technology.  This activity was part of a large
multi-disciplinary, multi-university project funded jointly by Ballard
and Mitacs (a Canadian research network dedicated to funding research
projects at the interface between mathematics and industry).  I joined
this project as a postdoctoral fellow in 1998 and continued as a faculty
investigator until roughly 2005.

The primary aim of the project was to develop mathematical models of fuel
cell processes and components that could support the Company's internal
R\&D efforts.  Space doesn't permit a detailed description of the wide
range of problems investigated and so I will only provide a 
high-level overview of a few main modelling efforts:
\begin{itemize}
\item A 2D model of heat and mass transport in a \mydefn{unit cell} that
  captures flow of a multi-component, multi-phase (gas-liquid) mixture
  through a porous fuel cell electrode.
\item Reduced models that focus solely on electrical coupling between
  unit cells and can be used to validate the results of cyclic
  voltammetry and other similar diagnostic tests.
\item Reduced models of an entire fuel cell that consist of algebraic
  equations, ODEs and simplified PDEs.  The goal here was to construct
  models that capture the essential physics but are simple enough to
  permit development of efficient algorithms for simulating fuel
  cell stacks (consisting of hundreds of unit cells connected in
  series).
\item A detailed model of the membrane that lies at the heart of the PEM 
  fuel cell, and whose complicated physics is still not very well
  understood.
\item Nano-scale models for the catalyst layer, which is a complex
  multi-scale composite material consisting of a porous
  electrically-conducting electrode, impregnated by the ionomer membrane
  and platinum catalyst.
\end{itemize}
These and other projects led to a host of academic publications (a few
of which I was personally involved with) that are detailed in a
comprehensive review article on the mathematics of PEM fuel cells
co-authored by the project's principal
investigators~\cite{promislow-wetton-2009}.

Before closing, I want to draw attention to the issue of
\emph{intellectual property} (or IP) which is an important aspect of
many industrial collaborations.  Ballard's business depends critically
on their PEM fuel cell-related inventions and patents, and so they are
naturally extremely sensitive to IP.  It was initially a struggle to
craft an IP agreement that was acceptable to the Company and the
universities involved, and the ultimate solution was both elegant and
surprisingly simple: the mathematics and the algorithms belonged to us
mathematicians, whereas the data and any simulations run on the data
were owned by the Company.  This left us relatively free to publish
mathematical results in the open literature with only a short delay for
the Company to review drafts of papers.  I have come to believe that
mathematics has a significant advantage over other disciplines (such as
engineering) in the sense that \emph{mathematical IP} is often
considered by companies to be of less obvious concern \dots\ although
perhaps it should be!

\leavethisout{
\subsection{Traffic flow: A serendipitous conversation}
\label{sec:traffic}

THIS IS MIbI!  In 2009, a previous graduate student contacted 

Traffic flow (2009): conversation with a company leading to two MSc
thesis projects but no collaboration (as
yet)~\cite{bowlby-2011, wiens-2011, wiens-stockie-williams-2013}.
}

\subsection{Atmospheric pollutant dispersion: Graduate student
  internships}
\label{sec:atmos}

I will next describe a collaboration with Teck, which is a
multi-national mining corporation headquartered in Vancouver.  This
project was initiated in 2005 by a \myquote{cold call} to SFU from an
Environmental Superintendent at Teck's lead-zinc smelter in Trail,
BC\footnote{In my experience, such cold calls (in either direction) are
  very rarely successful, which makes the manner in which this
  collaboration was initiated particularly notable.}.  The Company had
taken a series of ground-level measurements of particulate material
accumulated in \emph{dust-fall jars} around the site, and wanted to know
how they could exploit this data to estimate the corresponding airborne
contaminant emissions from a certain set of sources that they were
unable to measure directly.  The Company's goal was to provide
additional backing for (and perhaps even improvements to) the
engineering estimates of emissions that Teck is legally obliged to
report annually to Environment Canada.

\begin{figure}[hbt]
  \small
  \centering
  \begin{tabular}{ccccc}
    \includegraphics[width=0.29\textwidth]{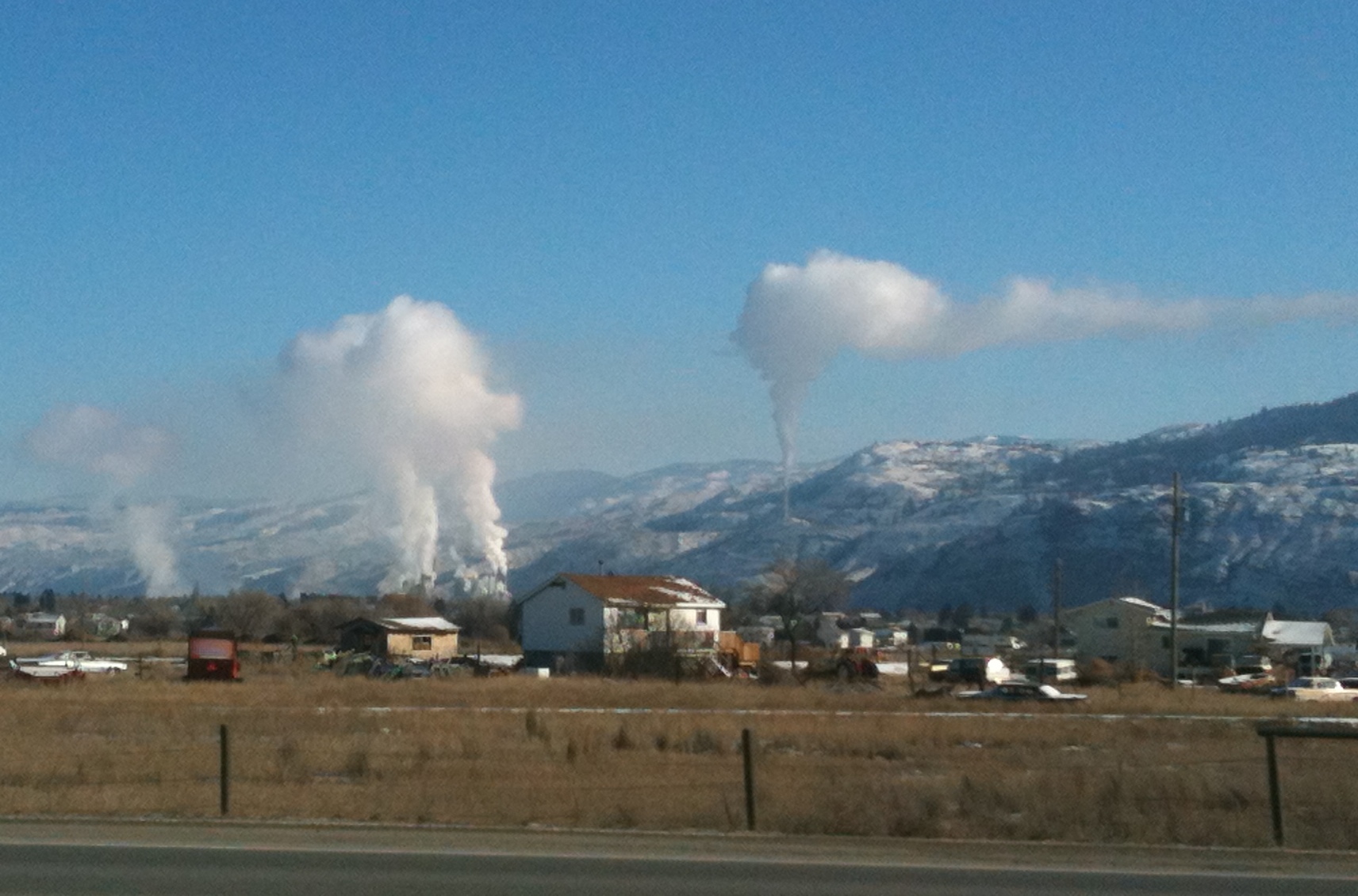}
    & \qquad & 
    \includegraphics[width=0.30\textwidth]{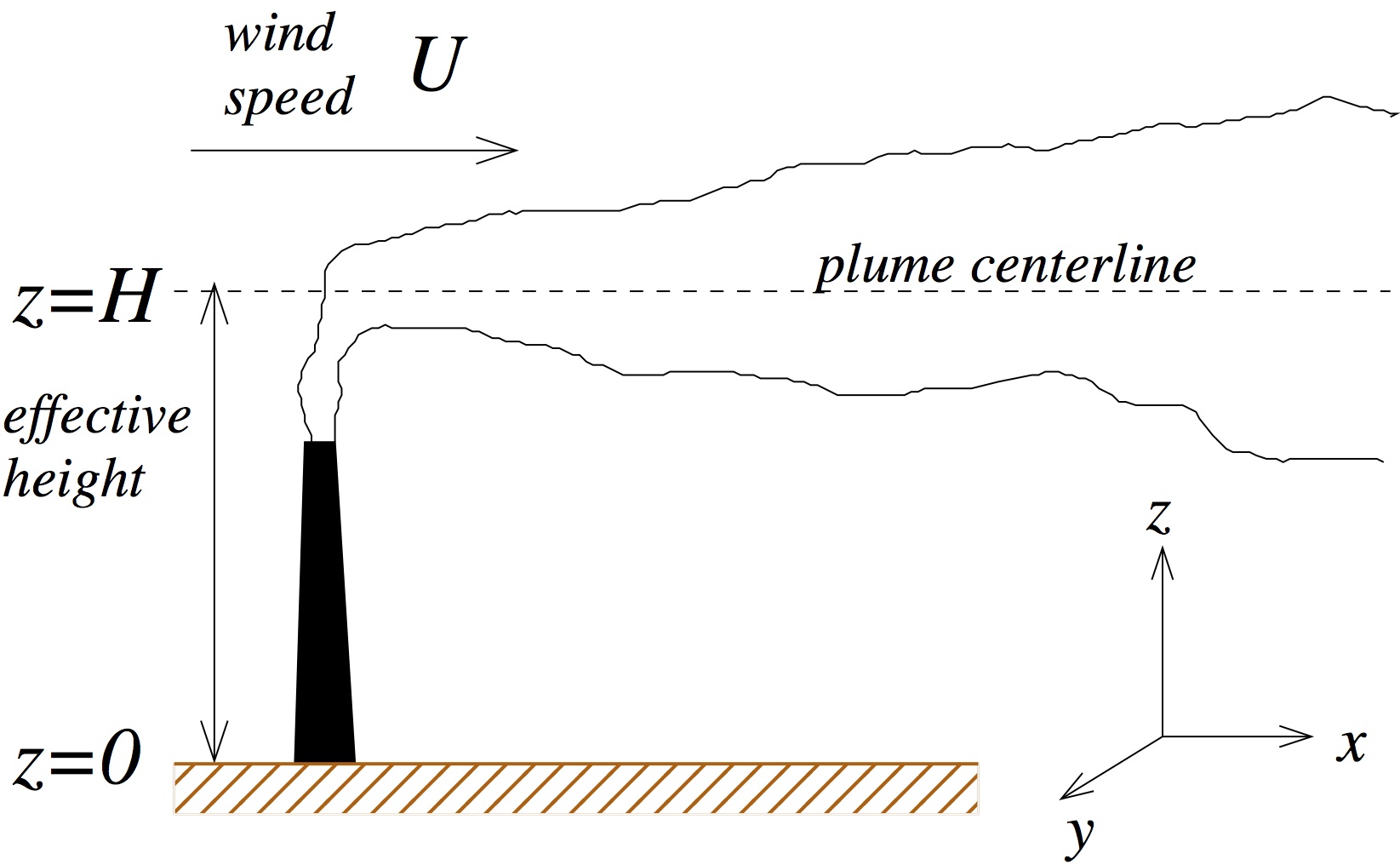}
    & \qquad & 
    \includegraphics[width=0.27\textwidth]{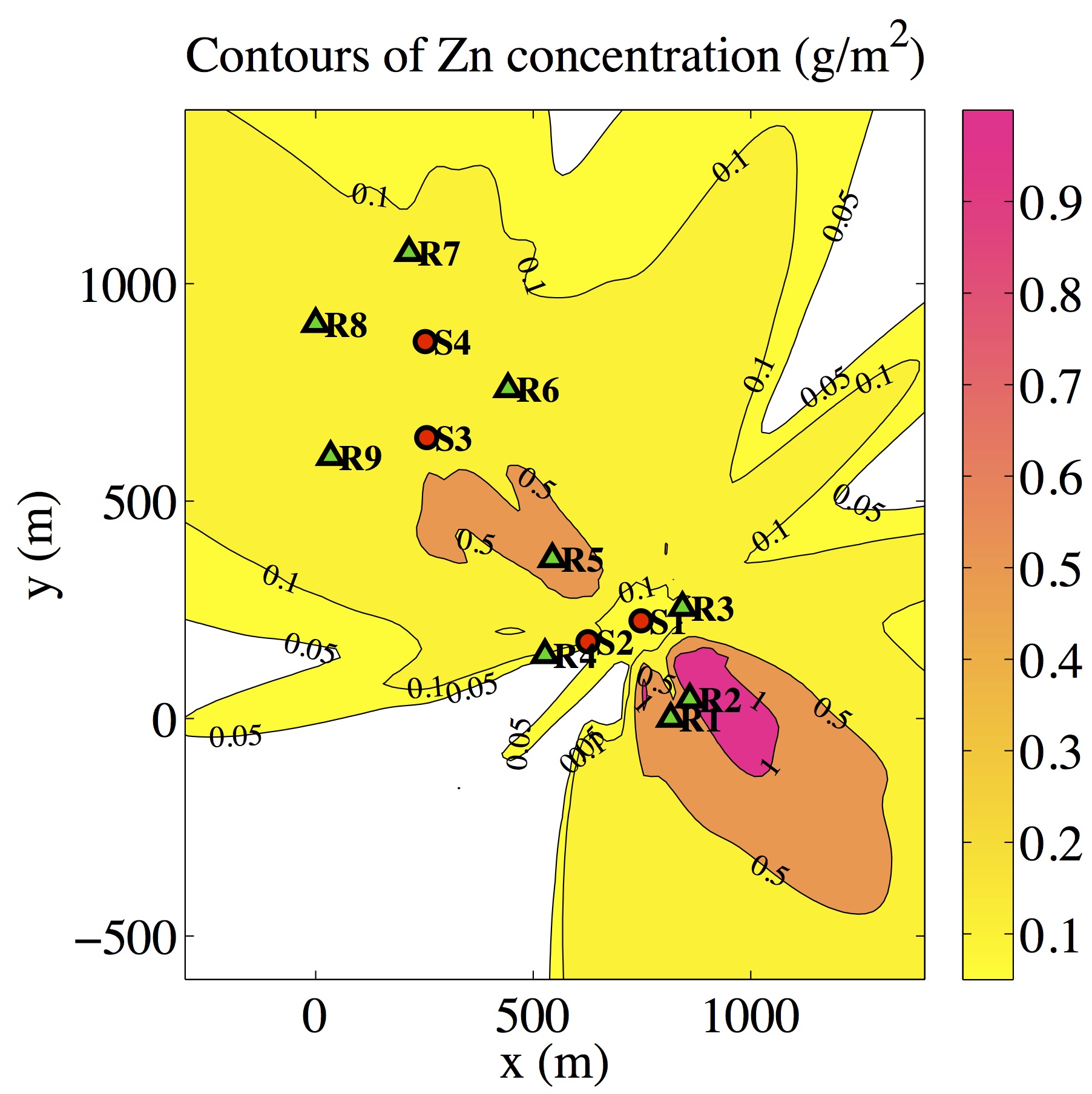}
    \\
    (a) & & (b) & & (c)
  \end{tabular}
  \caption{(a) Plumes generated by a pulp mill in Kamloops, BC.  (b)
    Schematic of a contaminant plume emitted from a source at location
    $(0,0,H)$ in a steady wind with velocity $(U,0,0)$. (c) Contours of
    ground-level zinc concentration over the smelter site, with four
    sources (S$n$, red circles) and nine dust-fall measurements (R$n$,
    green triangles).}
  \label{fig:atmos}
\end{figure}

Inspired by the fledgling Mitacs Graduate Internship Program (now known
as \emph{Mitacs Accelerate}) I convinced Teck to co-fund a Master's
student intern (E.~Lushi).  Within a few short months, she came up with a
surprisingly simple approach based on the \emph{Gaussian plume solution}
for contaminant concentration
\begin{gather*}
  C(x,y,z) = \frac{Q}{2\pi U \sigma_y \sigma_z}
  \exp\left(-\frac{y^2}{2\sigma_y^2} \right)
  \; \Bigg[ \exp\left( -\frac{(z-H)^2}{2\sigma_z^2} \right)
  + \exp\left(-\frac{(z+H)^2}{2\sigma_z^2} \right) \Bigg],
\end{gather*}
which is incidentally a very nice illustration of the application of
Green's functions or Laplace transform methods to solving the
advection-diffusion equation.  Referring to the problem geometry
depicted in Figure~\ref{fig:atmos}(a), the symbols appearing in this
equation are the emission rate $Q$, source height $H$, constant wind
velocity $U$ (in the $x$-direction), and dispersion coefficients
$\sigma_{y,z}$.  A sample ground-level zinc concentration map is shown
in Figure~\ref{fig:atmos}(b).  Taking advantage of the key fact that
concentration (and hence also deposition) is linear in $Q$, multiple
sources may be linearly superimposed and integrated in time to yield
cumulative depositions that afford a direct comparison with dust-fall
measurements.  The original question of determining emission rates based
on deposition measurements is an inverse problem that can therefore be
reformulated as a simple application of linear least squares.

This work has continued through four follow-up internships with Teck,
giving rise to a number of publications~\cite{hosseini-nigam-stockie-2015,
  lushi-stockie-2010, stockie-2011} and a Master's
thesis~\cite{hosseini-mscthesis-2013} in which a finite volume
discretization of the advection-diffusion equation was used to validate
the results from the inverse problem.  Furthermore, a PhD student
(B.~Hosseini) is currently studying more advanced Bayesian techniques
for source inversion problems and their application to both atmospheric
and groundwater contamination on the Teck smelter site.

\subsection{Maple sap exudation: A sweet project in search of a partner}
\label{sec:maple}

This last case study has an interesting background story that is worth
relating.  My first exposure to the science behind maple sap was
stumbling across an article in the tree physiology literature that
described a long-standing controversy over the bio-physical mechanisms
driving sap exudation during the spring thaw in trees like
maple\footnote{\mydefn{Exudation} refers to the unusual ability of
  certain deciduous trees such as maple, birch and walnut to generate
  stem pressure in a leafless state (i.e., in winter or early spring
  when transpiration is inactive).  When a sugar maple tree is
  \mydefn{tapped} during the harvest season, the exudation pressure is
  sufficient to cause the sweet sap to seep out in large enough
  quantities that it can be harvested commercially.}.  Having spent my
childhood in prime sugarbush territory in southern Ontario, I was
immediately intrigued by this quintessentially Canadian problem begging
for a mathematical solution.  However, without any existing connections
to tree physiologists or the maple syrup industry, I took the problem no
further (although as an inveterate problem-hunter I naturally filed it
away for future reference).  As luck would have it, a few years later I
encountered a call for proposals to the Research Fund of the North
American Maple Syrup Council -- a non-profit network of local industry
associations representing mostly small maple syrup producers.  After
submitting a proposal for a modest research grant in 2011, and
leveraging funds from a Mitacs postdoctoral fellowship program, I was
thrilled (indeed, even a bit surprised) when both proposals were successful.
The mathematical study of maple sap was on!

The first phase of the project involved developing a microscale model
for the \emph{Milburn-O'Malley process}, which is a freeze-thaw
mechanism that operates at the scale of individual wood cells and is
currently the most favoured hypothesis to account for maple's ability to
generate positive stem pressure.  This process involves an interplay
between liquid/frozen sap and trapped gas within two types of wood cell
called fibers and vessels (refer to Figure~\ref{fig:fibves}):
\begin{figure}[hbt]
  \small
  \centering
  \begin{tabular}{ccc}
    \includegraphics[width=0.2\textwidth]{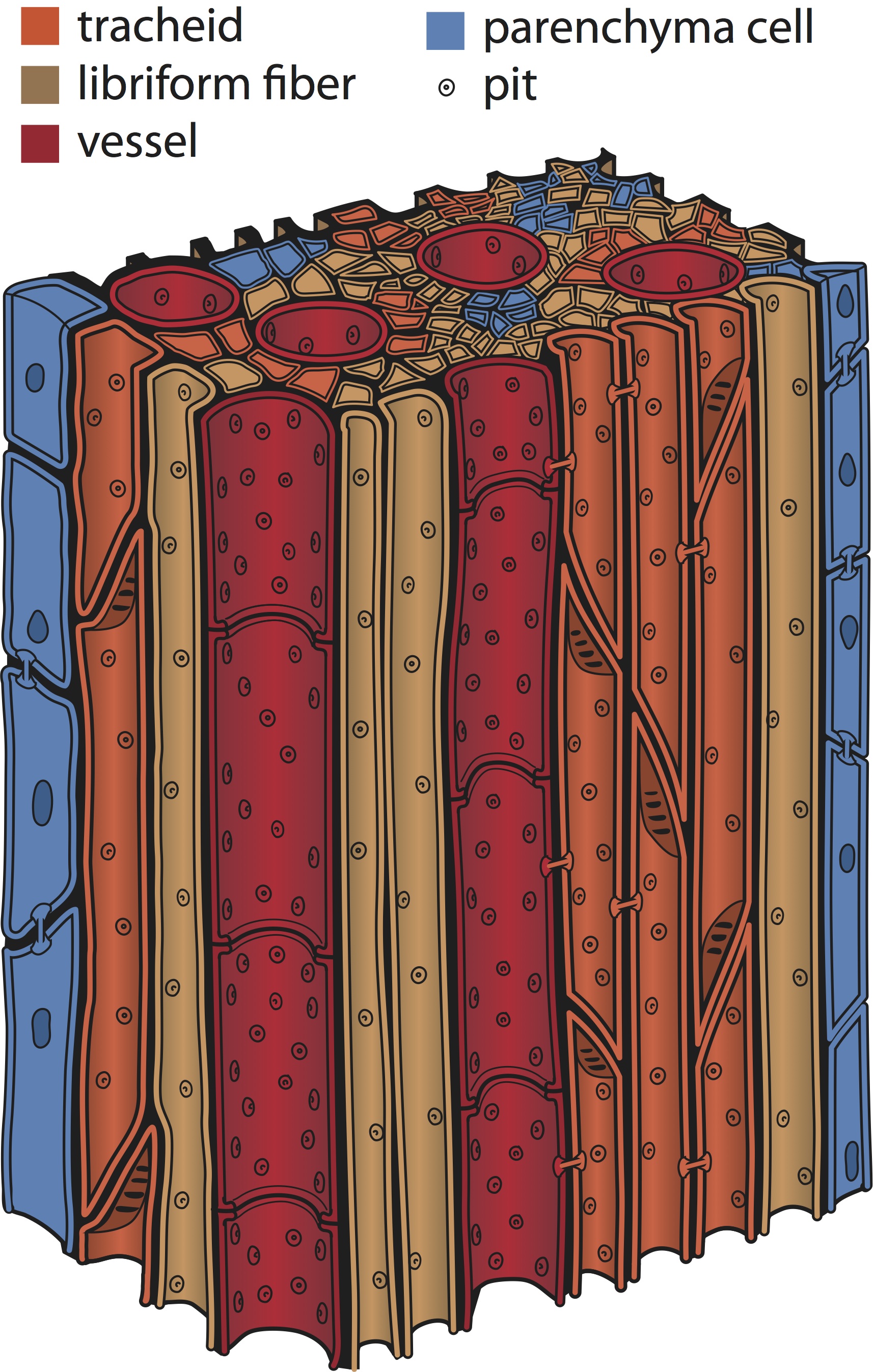}
    & \qquad\qquad & 
    \includegraphics[width=0.34\textwidth]{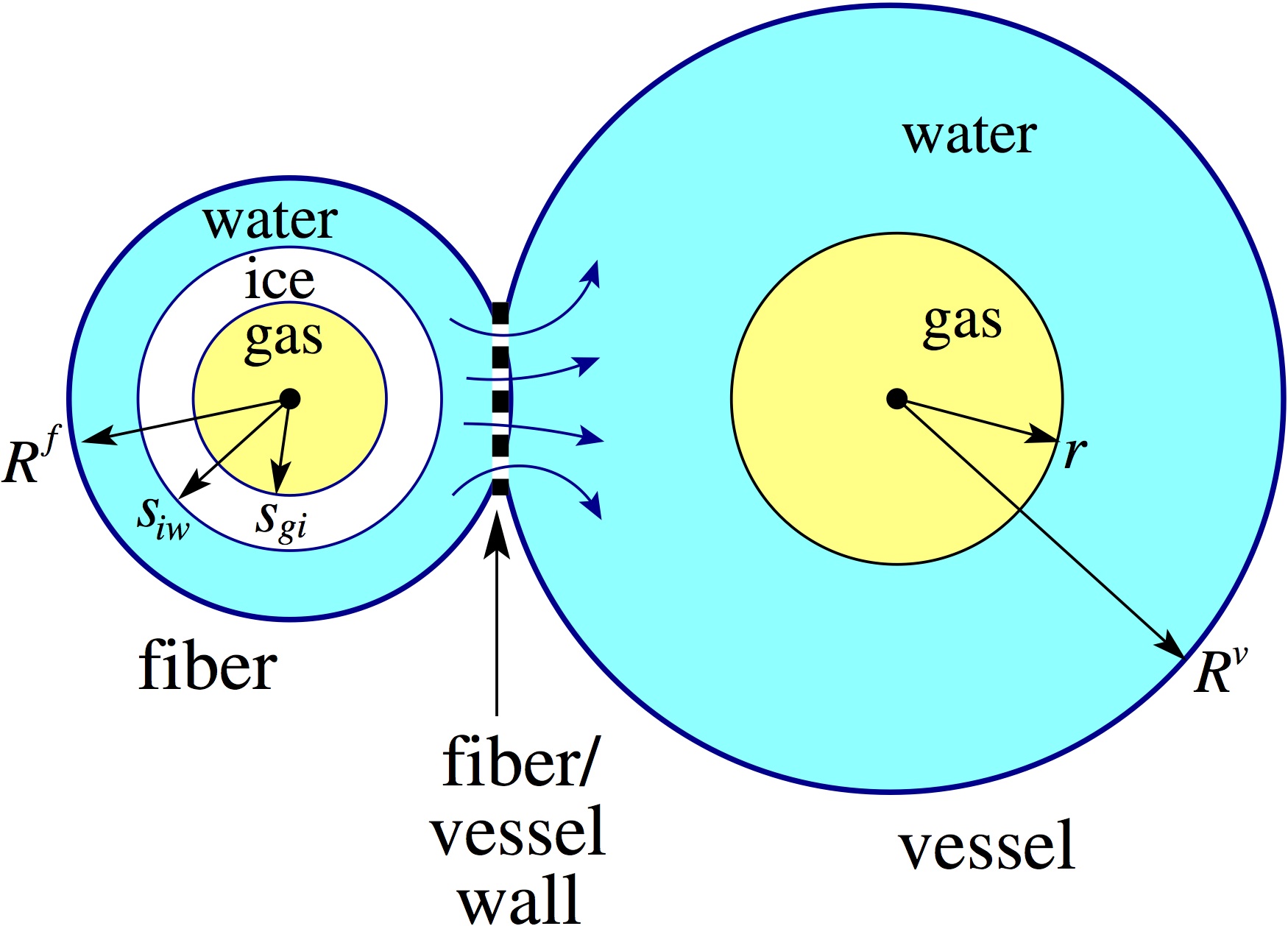}\\
    (a) & & (b)
  \end{tabular}
  \caption{%
    (a) 3D view of hardwood structure, showing vessels surrounded by
    (libriform) fibers and tracheids.
    %
    %
    (b) 2D cross-section through a fiber--vessel pair during the thawing
    process, showing the liquid, gas and ice regions, along with the
    phase interfaces.}
  \label{fig:fibves}  
\end{figure}
\begin{itemize}
\item In late autumn as temperatures drop below freezing, ice crystals
  form on the inner wall of the mostly gas-filled fiber cells, which
  compresses the gas trapped inside.
\item During the sap harvest season when temperatures rise above
  freezing, the ice layer melts and the pressure of the trapped gas
  drives the melted sap through the fiber walls into the vessels, hence
  causing the elevated pressures observed during sap exudation.
\item Pressures are sufficiently high that significant quantities of gas 
  dissolve within the sap.
\item Flow between the fiber and vessel passes through a selectively
  permeable membrane that sustains an osmotic pressure difference.
\end{itemize}
Along with postdoctoral fellow M.~Ceseri, I developed a mathematical
model of the thawing half of the process that consists of a coupled
system of nonlinear differential-algebraic equations incorporating heat
diffusion, porous media flow, Stefan conditions for phase interface
motion, osmosis, gas dissolution, and mass conservation.  By employing a
quasi-steady approximation for temperature (along with other simplifying
assumptions) we obtained a much simpler system of three ODEs governing
the location of the various gas-liquid
interfaces~\cite{ceseri-stockie-2013} that captures the expected
exudation behaviour.  Work is currently underway (with a second postdoc,
I.~Graf) to develop a corresponding model for the freezing half of the
process.

This project is a superb example of how an industrial collaboration can
be a rich source of stimulating mathematical problems.  For example, the
phase interface motion in our fiber-vessel model exhibits a non-trivial
separation of time scales that can be captured very accurately using an
asymptotic analysis~\cite{ceseri-stockie-2014}.  Furthermore, the
multi-scale wood cell structure lends itself naturally to methods of
periodic homogenization~\cite{graf-ceseri-stockie-2015}.  These and
other mathematical analyses support our two ultimate goals of (1)
addressing outstanding fundamental questions related to the biophysical
mechanisms governing sap exudation, and (2) applying these insights to
problems more relevant to the maple syrup industry such as developing
improved sap harvesting strategies or explaining the impact of climate
change on sap yields.  There is certainly no shortage of mathematical
problems related to maple sap that will keep the author and his
students/postdocs busy for years to come.

\leavethisout{
\emph{
  \begin{itemize}
  \item Describe interactions with farmers at NAMSC meetings.  
  \item Mention agriculture.
  \item Ceseri currently employed with Sportello Matematico
    per l'Industria Italiana (SMII), an organization that aims to build
    links between mathematicians and industry in Italy.
  \end{itemize}
}
}


\section{Conclusions}
\label{sec:conclude}

I hope that the preceding examples have managed to convey some of my own
enthusiasm for industrial problems, as well as convincing you of the
opportunities available to mathematicians within the field of research I
have termed \mydefn{mathematics for industry}.  Industrial
collaborations can be a rich source of interesting and challenging
mathematical research problems that lead to many spin-offs of an
academic and pedagogical nature.

It is worth highlighting a few recurring themes in this article that
lead naturally to some words of advice for academics who might be
interested in undertaking work at the interface between mathematics and
industry:
\begin{itemize}
\item Interesting and novel mathematical problems can arise in 
  unexpected places, and the people best suited to recognize and exploit
  these opportunities are mathematicians who have an interdisiplinary
  training and a broad knowledge of mathematics.

\item Serendipity plays an important role in identifying fruitful
  industrial collaborations.  One must be vigilant for opportunities and
  also willing and able to communicate with (potential) partners in
  clear and non-technical language.

\item There is a huge demand (and an even larger need, that is
  frequently unrecognized) for advanced mathematics in industry. It is
  often difficult for companies to engage with academics to gain entry
  to this expertise, and so anything that we as academic mathematicians
  can do to lower these barriers to access can be of huge mutual
  benefit.  A few suggestions are to create a web page that explains
  your work in lay terms, to introduce yourself to non-academic
  participants at conferences, to attend trade shows or other events
  that attract industrialists, and to give public or outreach lectures.

\item Be flexible, agile and willing to explore new areas.  The time
  investment required to move into a new sub-field of mathematics or
  application area introduces some degree of risk and delay for an
  academic mathematician, but the potential rewards are substantial.

\item Industrial collaborations can provide excellent opportunities for
  students and postdocs in terms of novel research projects and
  professional development.  Problems arising from these collaborations
  are often suitable for incorporating into undergraduate classes or
  outreach activities as ``modules'' that poignantly illustrate the value
  of mathematics to industry.

  \leavethisout{
  \item There is a sizable and growing industrial mathematics community
    world-wide, as documented in the OECD Global Science Forum
    report~\cite{oecd-gsf-2008}.  Many opportunities exist for attending
    conferences (organized by societies like ECMI and SIAM), participating in
    study groups (see \url{http://miis.maths.ox.ac.uk/how}), and
    publishing in targeted journals (such as \emph{Journal of Mathematics
      in Industry} from Springer and \emph{MICS Journal:
      Mathematics-in-Industry Case Studies} from Fields Institute).
  }
\end{itemize}
Finally, let me close with a few words of inspiration regarding
industrial mathematics from John Ockendon, who expresses his
\myquote{confidence in the intellectual viability of an activity that
  always has and surely always will spring mathematical surprises at a
  rate that could never be matched by most academic mathematicians
  pursuing their trade in the traditional way}~\cite{ockendon-2008}.


\providecommand{\noopsort}[1]{}

\end{document}